\documentclass[journal]{IEEEtran}

\usepackage{amsthm,amssymb,amsmath}
\usepackage{graphicx,tikz,pgfplots}
\usepackage[hidelinks]{hyperref}
\usepackage{subfigure}
\usepackage{dsfont}
\usepackage{enumerate}
\usepackage{textcomp}
\usepackage{empheq}
\usepackage{comment}
\usepackage{cleveref}
\usepackage[short,nocomma]{optidef}
\usepackage{booktabs} 
\usepackage{cite}

\theoremstyle{plain}
\newtheorem{lemma}{Lemma}
\newtheorem{theorem}{Theorem}
\newtheorem{corollary}{Corollary}

\theoremstyle{definition}
\newtheorem{assumption}{Assumption}

\theoremstyle{remark}
\newtheorem{remark}{Remark}
\newtheorem*{example*}{Illustrative Example}

\DeclareMathOperator{\Tr}{\mathbf{Tr}}

\newcounter{mytempeqncnt}

\newcounter{storesubeqcounter}

\def\BibTeX{{\rm B\kern-.05em{\sc i\kern-.025em b}\kern-.08em
    T\kern-.1667em\lower.7ex\hbox{E}\kern-.125emX}}
    
\makeatletter
\renewcommand*\env@matrix[1][*\c@MaxMatrixCols c]{%
  \hskip -\arraycolsep
  \let\@ifnextchar\new@ifnextchar
  \array{#1}}
\makeatother

\begin{document}
\title{Direct data-driven control of LTV systems}

\author{B. Nortmann and T. Mylvaganam  
\thanks{B. Nortmann  and T. Mylvaganam are with the Department of Aeronautics, Imperial College London, London SW7 2AZ, UK.}
}

\maketitle

\begin{abstract}
Considering discrete-time linear time-varying systems with unknown dynamics, controllers guaranteeing bounded closed-loop trajectories, optimal performance and robustness to process and measurement noise are designed via convex feasibility and optimisation problems involving purely data-dependent linear matrix inequalities. 
For the special case of periodically time-varying systems, infinite-horizon guarantees are achieved based on finite-length data sequences. 
\end{abstract}


\section{Introduction}\label{sec:intro}

\textit{Direct data-driven control} methods, which aim to control a system directly using data, without explicitly identifying a system model, have recently attracted significant interest (see, \emph{e.g.} \cite{Hou2017}). 
A central question in direct data-driven control is how to substitute a system model with data. For linear time-invariant (LTI) systems, a recent line of research addresses this question via Willems et al.'s \textit{fundamental lemma} \cite{Willems}.
The result is used in \cite{Coulson}, \cite{Berberich2021} to replace the system model and initial conditions in the context of model predictive control with data. In \cite{DePersis2020} it is used to derive a data-driven representation of closed-loop systems under static state feedback, where the controller itself is parametrised using data only.
This can be used to formulate and solve the stabilisation \cite{DePersis2020}, linear quadratic regulator (LQR)  \cite{DePersis2020,Rotulo} and suboptimal control \cite{vanWaarde2020} problems in terms of data-dependent linear matrix inequalities (LMIs). 
Extensions have been proposed for data from multiple data sets \cite{vanWaardeMultiple}, certain classes of nonlinear systems \cite{Guo2020DatadrivenSO,Strsser2020DataDrivenSO,Alsalti2021DataBasedSA},
linear parameter-varying systems \cite{Verhoek2021FundamentalLF} and switched systems \cite{Rotulo2021OnlineLO}.
In \cite{Nortmann2020DataDrivenCO}, the direct data-driven control framework originally presented in \cite{DePersis2020} (for LTI systems) is extended to linear time-varying (LTV) systems. Time-varying systems arise in a variety of practical problems and LTV models emerge, for instance, when linearising nonlinear systems around a trajectory or time-varying operating point \cite{LTVbook}. 
The demand for model-free control approaches for LTV systems is apparent in the literature (see, \emph{e.g.} \cite{Fong2018, Liu2017, Pang2018, Baros2020}). 
In this paper, we provide a complete analysis of the preliminary results in \cite{Nortmann2020DataDrivenCO}. Moreover, we extend the results to LTV systems affected by both measurement and process noise and provide insights for the special case of periodically time-varying systems.
In contrast to the related results \cite{Verhoek2021FundamentalLF}, \cite{Rotulo2021OnlineLO}, the presented data-driven methods are applicable to (linear) \emph{arbitrarily} time-varying systems and do not rely on any assumptions or prior knowledge of the system structure or parameter variation. However, it is shown how such knowledge can be exploited for the special case of periodically time-varying systems.
Challenges associated with direct data-driven control in the presence of noise are addressed for certain classes of control problems involving LTI systems in \cite{DePersis2020, DePersis2021, Berberich2020, vanWaarde2021, Bisoffi2021}. Most works in this context focus on process noise only, or consider process noise and measurement noise separately.
The results presented herein are inspired by \cite{DePersis2020,Berberich2020}, and can be considered as an LTV equivalent. The main difference apart from the extension to LTV systems - which itself introduces new challenges and requires a different approach to parametrise unknown systems - is 
that we incorporate \emph{both} measurement and process noise in a single formulation and study the behaviour of the system in closed-loop with feedback on the noisy state measurements.

The paper is organised as follows. In Section~\ref{sec:prelim} some preliminaries are provided.
In Section~\ref{sec:no_noise} we consider noise-free LTV systems and show that state feedback control laws (guaranteeing a decreasing bound on the closed-loop trajectories or solving the time-varying LQR problem) can be designed via data-dependent semidefinite programmes (SDPs).
The problem of designing controllers with robustness guarantees directly using noisy data is addressed in Section~\ref{sec:H_inf}. 
In Section~\ref{sec:periodic}, we specialise the results to the class of periodically time-varying systems.
Concluding remarks are provided in Section~\ref{sec:conclusion}.
 
\smallskip
\noindent
\textbf{Notation.} 
The sets of real numbers, integers and natural numbers are denoted by $\mathbb{R}$,  $\mathbb{Z}$ and  $\mathbb{N}$, respectively.
The zero matrix of appropriate dimension is denoted by $0$ and the $n\times n$ identity matrix by $I_n$.
Given a square matrix $A$, $\Tr (A)$ denotes its trace, and $A \succ 0$ ($A \succeq 0$) denotes that $A$ is positive definite (positive semi-definite). In matrix inequalities $\star$ denotes blocks (or matrices), which can be inferred by symmetry. 
The block diagonal stacking of matrices $A$ and $B$ is written as $\text{diag} \left( A, B \right)$. 
Given a vector $v\in \mathbb{R}^n$, $\|v\|$ denotes its Euclidean norm and given a matrix $M\in\mathbb{R}^{m\times n}$, $\|M\|$ denotes the induced 2-norm of $M$. 
Given a signal $z: \mathbb{Z} \rightarrow \mathbb{R}^{\sigma}$ the sequence $ \{ z(k), \ldots, z(k+T) \} $ is denoted by $z_{ \left[ k,k+T \right]}$ with $k,T \in \mathbb{Z}$ and we denote $|z|_k = \text{sup} \left\lbrace\|z(j)\|, 0 \leq j \leq k \right\rbrace \leq \infty$. The space of square-summable sequences is denoted by $\ell_2$. A function $\gamma: \mathbb{R}_{\geq 0} \rightarrow \mathbb{R}_{\geq 0}$ is a class $\mathcal{K}$-function if it is continuous, strictly increasing and $\gamma(0) = 0$.

\section{Preliminaries and problem definition}\label{sec:prelim}
Consider a discrete-time LTV system, described by 
\begin{subequations}
\begin{equation}
x(k+1) = A(k)x(k) + B(k)u(k) + d(k), 
	\label{eq:noisy dynamics state}
\end{equation}
where $x \in \mathbb{R}^n$ is the state of the system, $u \in \mathbb{R}^m$ is the control input, $d \in \mathbb{R}^n$ denotes an unknown additive system noise and $A(k)$ and $B(k)$ denote the unknown time-varying dynamics and input matrices of appropriate dimensions, respectively.
Suppose that available state measurements $\zeta \in \mathbb{R}^n$  are corrupted by measurement noise $v \in \mathbb{R}^n$, \emph{i.e.} 
\begin{equation}
	\zeta(k) = x(k) + v(k).
	\label{eq:noisy measurements}
\end{equation}
\label{eq:noisy system dynamics}%
\end{subequations}
Our objective is to design controllers of the form
\begin{equation}
	u(k) = K(k) \zeta(k),
	\label{eq:noisy state-feeback}
\end{equation}
for the unknown LTV system \eqref{eq:noisy system dynamics} based solely on measurements of the (noisy) state and input of the system, such that certain guarantees hold for the resulting closed-loop system
\begin{equation}
x(k+1) = \bigg( A(k)+B(k)K(k) \bigg) x(k) +  B(k)K(k)v(k) + d(k).
\label{eq:noisy closed-loop}
\end{equation}
If $d(k) = 0$ and $\zeta(k) = x(k)$, for all $k$, \emph{i.e.} in the noise-free case (which is considered in Section~\ref{sec:no_noise}), the dynamics of the unknown LTV system \eqref{eq:noisy system dynamics} simplify to 
\begin{equation}
    x(k+1) = A(k) x(k) + B(k) u(k),
    \label{eq:systemdynmaics}
\end{equation}
and the closed-loop system under state feedback with $u(k) = K(k) x(k)$, is described by
\begin{equation}
    x(k+1) = \bigg( A(k) + B(k)K(k) \bigg) x(k).
    \label{eq:closed loop}
\end{equation}
We consider the direct data-driven solution of several classical control problems involving the unknown systems \eqref{eq:systemdynmaics} or \eqref{eq:noisy system dynamics}.
To this end, we summarise some preliminaries regarding data-driven control of (noise-free) LTV systems. These are extended to the general case including noise in Section~\ref{sec:H_inf}.
Throughout the paper, we 
consider the following standing 
assumption.
\begin{assumption}
It is possible to gather an ensemble of $L \in \mathbb{N}$ input-state data sequences capturing the \emph{same time-varying behaviour} of the unknown LTV system over $T+1$ time instances, with $T \in \mathbb{N}$, \emph{i.e.} if data sequence $j$ covers the time interval $k = k_j, \ldots, k_j + T$, for $j = 1, \ldots, L$, then for all $l=1, \ldots, L$,\footnote{In the remainder of the paper we refer to each interval capturing the time-variation of interest as $k = 0, \ldots,T$, \emph{i.e.} we define $k_j = 0$, for $j = 1, \ldots,L$.} $\{A(k_j), \ldots, A(k_j+T-1)\} = \{A(k_l), \ldots, A(k_l+T-1)\}$ and $\{B(k_j), \ldots, B(k_j+T-1)\} = \{B(k_l), \ldots, B(k_l+T-1)\}$. 
\label{as:conditions}
\end{assumption}
\begin{remark}
Assumption~\ref{as:conditions} is similar to requirements typically encountered in ensemble methods for LTV system identification (see \emph{e.g.} \cite{Verhaegen1995}). 
It is readily satisfied by systems arising in a variety of applications, including biomedical systems, nonlinear systems linearised along a trajectory and periodically varying systems (which are addressed in Section~\ref{sec:periodic}).  
Variations in environmental conditions that may result in different time-variations affecting each experiment in an ensemble can be considered as process noise, which is addressed in Section~\ref{sec:H_inf}. 
\label{re:context to Assumption 1}
\end{remark}
\noindent The $L$ input-state data sequences can be obtained via a sequence of physical experiments or via simulations\footnote{Herein, we refer to the act of data collection as ``experiment'', regardless of whether the data is collected via experiments or simulations.}. Considering the system \eqref{eq:systemdynmaics}, let $u_{d,j,\left[ 0, T-1 \right]}, \ x_{d,j,\left[ 0, T \right]},
$ represent input-state data collected during the $j^{\mathrm{th}}$ experiment, for $j = 1, \ldots, L$. 
While the specific experiment is indicated by the subscript $j$, the subscript $d$ highlights that the input-state sequences contain \emph{measured data} samples. 
Consider the matrices
\vspace{-0.5cm}
\begin{subequations}
\begin{equation}
    \begin{split}
        X(k) & = \left[ x_{d,1}(k), x_{d,2}(k), \ldots, x_{d,L}(k) \right], \\
\end{split}
\label{eq:X(k)}
\end{equation}
for $k = 0, \ldots T$, and
\begin{equation}
    \begin{split}
        U(k) & = \left[ u_{d,1}(k), u_{d,2}(k), \ldots, u_{d,L}(k) \right],
    \end{split}
\label{eq:U(k)}
\end{equation}
\end{subequations}
for $k = 0, \ldots T-1$, which combine the data from all $L$ experiments at each time step. Note that the data matrices $X(k)$ and $U(k)$ satisfy 
\begin{equation}
        X(k+1)  = A(k) X(k) + B(k) U(k) 
         = \left[ A(k) \  B(k) \right] 
        \begin{bmatrix}
            X(k) \\
            U(k)
        \end{bmatrix},
    \label{eq: data matrix dynamics}
\end{equation}
for $k=0,\ldots,T-1$.
Suppose the rank condition\footnote{The condition \eqref{eq:rank condition} can always be verified from the measured data. A necessary condition for \eqref{eq:rank condition} to hold is that $L \geq n+m$.} 
\begin{equation}
    \mathrm{rank} 
    \begin{bmatrix}
        X(k) \\
        U(k)
    \end{bmatrix} = n + m, 
    \label{eq:rank condition}
\end{equation}
holds for $k = 0,\ldots,T-1$. Then, the closed-loop system \eqref{eq:closed loop} can equivalently be represented as
\begin{align}
    x(k+1) & = X(k+1)G(k)x(k), \label{eq: param closed loop} 
\end{align}
where $G(k) \in \mathbb{R}^{L \times n}$ satisfies
\begin{equation}
    \begin{bmatrix}
        I_n \\
        K(k)
    \end{bmatrix} = 
    \begin{bmatrix}
        X(k) \\
        U(k)
    \end{bmatrix} G(k),
    \label{eq: G definition}
\end{equation}
for $k = 0, \ldots, T-1$, see \cite{Nortmann2020DataDrivenCO}.
Hence, the sequence of control gains $K(k)$ is parametrised using data through the identity \eqref{eq: G definition}.
Thus, the matrices $G(k)$, for $k=0,\ldots,T-1$, can be seen as decision variables, which can be used for identification-free design of state feedback controllers. 

\begin{remark}
To utilise the data-driven system representation \eqref{eq: param closed loop}, \eqref{eq: G definition},
it is required that 
\eqref{eq:rank condition} holds for all $k = 0, \ldots, T-1$.
Thus, each input-state data sequence $j$, $j=1,\ldots,L$, in the ensemble must start at different initial conditions $x_{d,j}(0)$.
If this is infeasible, a common starting point can be considered as the state at $k=-1$ and different inputs can be applied for each experiment to obtain different state data at $k=0$.
\label{re:different initial conditions}
\end{remark}

\section{Data-driven control of LTV systems - the Noise-free case}\label{sec:no_noise}
In this section, we utilise the direct data-driven system representation \eqref{eq: param closed loop}, \eqref{eq: G definition} to design feedback controllers for the unknown (noise-free) LTV system \eqref{eq:systemdynmaics} via the solution of convex optimisation problems involving LMI constraints.

\subsection{Bounded closed-loop trajectories} 
\label{subsec:bound}
Consider the problem of controlling the LTV system \eqref{eq:systemdynmaics} over a finite time horizon, with the aim of ensuring that the closed-loop trajectories remain close to the equilibrium throughout the considered horizon.
A solution to this problem is provided in the following statement.

\begin{theorem}
    Consider the system \eqref{eq:systemdynmaics} and suppose an ensemble of input-state data is available to form the matrices \eqref{eq:X(k)}, \eqref{eq:U(k)}, such that the rank condition \eqref{eq:rank condition} holds, for $k=0,\ldots,T~-~1$.
    Any sequences of matrices $Y(k)$, $P(k)=P(k)^\top$ satisfying
    \begin{subequations}
    	\begin{equation}
        \begin{bmatrix}
            P(k+1)-I_n & X(k+1) Y(k) \\
            Y(k)^\top X(k+1)^\top & P(k) 
        \end{bmatrix} \succeq 0, \label{eq:param stability inequality} \\ 
        \end{equation}
        \begin{equation}
        	X(k)Y(k) = P(k), \label{eq:param stability equality} \\
        \end{equation} 
        for $k = 0, \ldots, T-1$, and
        \begin{equation}
        	\eta I_n \preceq P(k) \preceq \rho I_n, \label{eq:P boundedness}
		\end{equation}          
		for $k = 0, \ldots, T$, where $\eta \geq 1$ and $\rho > \eta$ are finite constants, are such that the trajectories of the system \eqref{eq:closed loop}, with
    \label{eq:param stability condition}
    \end{subequations} 
    \begin{equation}
        K(k) = U(k)Y(k)P(k)^{-1},
        \label{eq:param control gain}
    \end{equation}
    for $k = 0, \ldots, T-1$, satisfy
    \begin{equation}\label{eq:traj_bound}
        \|x(k)\| \le \sqrt{\frac{\rho}{\eta}}\displaystyle\left(1-\frac{1}{\rho}\displaystyle\right)^{\frac{k}{2}}\|x(0)\|\,,
    \end{equation}
    for $k=0,\dots,T$.
    \label{th: param boundedness}
\end{theorem}
\begin{proof}
    To demonstrate the claim it is useful to consider the adjoint equation (see, \emph{e.g.} \cite[Chapter 3.1]{Ichikawa2001}) of the closed-loop system \eqref{eq:closed loop}. Namely, consider 
    \begin{equation}\label{eq:adjoint} 
        \xi(j) = A_{cl}(j)^\top \xi(j+1)\,, 
    \end{equation}
    and note that the solution to \eqref{eq:adjoint} starting from $\xi(k)$ is
    \begin{equation}\label{eq:adjoint_transit}
        \xi(j)= S_t(k,j)^\top\xi(k)\,, 
    \end{equation}
     for $j\le k$, where 
     $$
     S_t(k,j) = \begin{cases}
    A_{cl}(k-1)A_{cl}(k-2)\dots A_{cl}(j) \,, \quad &\text{for $j<k$} \,,\\ 
    I_n \,, \quad &\text{for $j=k$}\,,
    \end{cases}
     $$ 
     denotes the state transition matrix corresponding to the closed-loop system \eqref{eq:closed loop}, and where $A_{cl}(k) = A(k)+B(k)K(k)$. Let $\xi(k) \ne 0$ and suppose we can determine a sequence of matrices $P(k)$ satisfying the condition \eqref{eq:P boundedness} and 
     	\begin{equation}
              A_{cl}(k)  P(k) A_{cl}(k)^\top 
            -P(k+1)+I_n \preceq 0,
            \label{eq:stability condition}
        \end{equation}
    for $k=0, \ldots, T-1$. 
    Consider the quadratic function $\bar{V}_j := \bar{V}(j,\xi(j)) = \xi(j)^\top P(j) \xi(j)$, for $j=0, \dots, T$.
    It follows from  \eqref{eq:adjoint}, \eqref{eq:stability condition} and \eqref{eq:P boundedness} that $\bar{V}_{j+1} - \bar{V}_j \geq \|\xi(j+1)\|^2 \geq \frac{1}{\rho}\xi(j+1)^\top P(j+1) \xi(j+1),$ for $j=0, \dots, T-1$.
    Thus, we have that $\xi(j)^\top P(j) \xi(j) \leq  \left( 1-\frac{1}{\rho} \right)^{k-j}\xi(k)^\top P(k) \xi(k)$, for $j=0, \dots, T$, $j\le k \le T$.
    It then follows from \eqref{eq:P boundedness} that 
    $
    \eta \|\xi(j)\|^2 \leq \rho \left(1-\frac{1}{\rho}\right)^{k-j}\|\xi(k)\|^2\,, 
    $
    which, using \eqref{eq:adjoint_transit}, in turn yields 
    $
    \|S_t(k,j)^\top \xi(k) \|^2   \leq \frac{\rho}{\eta} \left(1-\frac{1}{\rho}\right)^{k-j} \|\xi(k)\|^2 \,, 
    $
    and 
    \begin{equation}\label{eq:stab} 
    \begin{split} 
    \|S_t(k,j)^\top &\|^2  = \|S_t(k,j)  \|^2 \leq \frac{\rho}{\eta} \left(1-\frac{1}{\rho}\right)^{k-j}\,, 
    \end{split} 
    \end{equation}  
    for $j=0, \dots, T$,
    $j\le k \le T$. 
    Noting that $x(k)~=~S_t(k,j) x(j)$, for $k \ge j$, \eqref{eq:stab} implies 
    \begin{equation} \begin{split} 
    \|x(k)\|^2 &= \|S_t(k,j)x(j)\|^2 \leq \|S_t(k,j)\|^2\|x(j)\|^2 \\& \leq \frac{\rho}{\eta}\left(1-\frac{1}{\rho}\right)^{k-j} \|x(j)\|^2 \,, 
    \end{split} 
    \end{equation}
     for $j=0,\dots,T$, $j\le k \le T$. Letting $j=0$, this yields \eqref{eq:traj_bound}. 
    Finally, using \eqref{eq: param closed loop}, \eqref{eq: G definition}, defining $Y(k):=G(k)P(k)$, and via the Schur complement, \eqref{eq:param stability inequality} is equivalent to \eqref{eq:stability condition}, if \eqref{eq:param stability equality} holds and the control gain is chosen as \eqref{eq:param control gain}.
\end{proof}

\subsection{Optimal control}\label{sec:LQR}
Consider system \eqref{eq:systemdynmaics} and the problem of finding the optimal control sequence $\left\{ u^*(0),u^*(1),\ldots, u^*(N-1) \right\}$ as a function of the state, which minimises the quadratic cost functional
\begin{multline}
    J \left( x(0), u(\cdot) \right) = x(N)^\top Q_f x(N) \\
    + \sum_{k=0}^{N-1} \left( x(k)^\top Q(k) x(k) + u(k)^\top R(k) u(k) \right),
    \label{eq: LQR cost}
\end{multline}
over the time horizon $N \in \mathbb{N}$, starting from the initial condition $x(0) = x_0$, with $Q_f = Q_f^\top \succeq 0$, $Q(k) = Q(k)^\top \succeq 0$ and $R(k) = R(k)^\top \succ 0$, for $k = 0,\ldots,N-1$. 
In \cite{Nortmann2020DataDrivenCO} it has been shown that this finite-horizon LQR problem can equivalently by solved via a convex programme. In the following statement, we combine this with the data-driven system representation \eqref{eq: param closed loop}, \eqref{eq: G definition} to formulate the time-varying LQR problem as a data-dependent SDP. 
\begin{theorem} \label{th:data LQR}
	Consider the system \eqref{eq:systemdynmaics} and suppose an ensemble of input-state data is available to form the matrices \eqref{eq:X(k)}, \eqref{eq:U(k)}, such that the rank condition \eqref{eq:rank condition} holds, for $k = 0, \ldots, N-1$.
	The optimal state feedback control gain sequence $\left\{ K^*(0), K^*(1), \ldots, K^*(N-1) \right\}$ solving the finite-horizon LQR problem with $u^\star(k) = K^*(k)x(k)$ is given by 
    \begin{equation}
        K^*(k) = U(k)H^*(k)S^*(k)^{-1},
        \label{eq:param gain}
    \end{equation}
    for $k = 0, \ldots, N-1$, with $H^*(k)$ and $S^*(k)$ the solution of
    \begin{mini!}[2]<b>
    	  {\mathcal{S,H,O}}{ \Tr \left( Q_f S(N) \right) \nonumber}{\label{eq:data-dependent LQR problem}}{}
    	  \breakObjective{+ \sum_{k=0}^{N-1} \bigg( \Tr \left( Q(k) S(k) \right) + \Tr \left(O(k) \right) 			\bigg) \label{eq:data-dependent LQR problem objective}}
    	  \addConstraint{S(0)}{\succeq I_n, \label{eq:data-dependent LQR problem inequality 1}}{ }
    	  \addConstraint{\begin{bmatrix} S(k+1)-I_n & X(k+1)H(k) \\ H(k)^\top X(k+1)^\top & S(k) \end{bmatrix}}			{\succeq 0, \label{eq:data-dependent LQR problem inequality 2}}
    	  \addConstraint{\begin{bmatrix} O(k) & R(k)^{1/2} U(k)H(k) \\ H(k)^\top U(k)^\top R(k)^{1/2} & S(k) 				\end{bmatrix}}{\succeq 0, \label{eq:data-dependent LQR problem inequality 3}}
    	  \addConstraint{S(k)}{=X(k)H(k), \label{eq:data-dependent LQR problem equality}}
    \end{mini!}
    for $k = 0,\ldots, N-1$, where $ 
            \mathcal{S} = \left\{ S(1), \ldots, S(N) \right\},$  
            $ \mathcal{H} = \left\{ H(0), \ldots, H(N-1) \right\}$ and
            $\mathcal{O} = \left\{ O(0), \ldots, O(N-1) \right\}$. 
\end{theorem}
\begin{proof}
    The proof lies in demonstrating that \eqref{eq:data-dependent LQR problem} is equivalent to the model-based convex programme corresponding to equation $(15)$ in \cite{Nortmann2020DataDrivenCO}.
    This follows by introducing \eqref{eq: param closed loop}, \eqref{eq: G definition} to the constraints, letting $H(k) := G(k)S(k)$, and taking the Schur complement of the nonlinear inequality constraints. 
\end{proof}

\section{Data-driven control of LTV systems - Robustness to noise} \label{sec:H_inf}
In practice, both the measurement and/or the system dynamics may be subject to noise.
In this section, we consider the problem of designing feedback controllers for the general (unknown) discrete-time LTV system \eqref{eq:noisy system dynamics}. 
To this end, we start by deriving a data-driven system representation of the form \eqref{eq: param closed loop}, \eqref{eq: G definition} using noise corrupted data. Namely, let $u_{d,j,\left[ 0, T-1 \right]}, \ \zeta_{d,j,\left[ 0, T \right]},$ denote input-output data collected during the $j^\text{th}$ experiment, for $j = 1,\ldots,L$.  
The data is arranged to form the matrices
\begin{equation}
        Z(k) = \left[ \zeta_{d,1}(k), \zeta_{d,2}(k), \ldots, \zeta_{d,L}(k) \right],
        \label{eq:Z(k)}
\end{equation}
for $k = 0, \ldots T$, representing the ensemble of noisy state measurements, and \eqref{eq:U(k)}, for $k = 0, \ldots T-1$, representing the ensemble of input measurements. Consider also a similar stacking of the corresponding samples of the system noise $d_d$, the state $x_d$ and the measurement noise $v_d$ (all of which are not measured) associated with the ensemble of experiments, represented by the matrices $D(k)$, for $k = 0,\ldots, T-1$, and \eqref{eq:X(k)} and $V(k)$, for $k = 0,\ldots, T$, respectively. Note that \begin{equation*}
\begin{split}
X(k+1) &= A(k)X(k) + B(k)U(k) + D(k), \\
Z(k) &= X(k) + V(k).
\end{split}
\end{equation*}
Suppose the rank condition
\begin{equation}
 \text{rank} \begin{bmatrix}
 Z(K) \\
 U(k)
 \end{bmatrix} = n+m,
 \label{eq:noisy rank condition}
 \end{equation}
holds for $k = 0,\ldots, T-1$. Then, the dynamics matrix of the closed-loop system \eqref{eq:noisy closed-loop} can equivalently be represented as
\begin{equation}
A(k) + B(k)K(k) = \bigg(Z(k+1)+W(k)\bigg) G(k),
\label{eq:noisy param closed loop}
\end{equation}
where $G(k)$ satisfies
\begin{equation}
\begin{bmatrix}
I_n \\
K(k) 
\end{bmatrix} = \begin{bmatrix}
Z(k) \\
U(k)
\end{bmatrix} G(k),
\label{eq:noisy G definition}
\end{equation}
for $k = 0,\ldots, T-1$, with
\begin{equation}
W(k) = A(k)V(k) - V(k+1) - D(k).
\label{eq:noisy R(k) definition}
\end{equation}
Assuming the unknown ensemble matrix $W(k)$, containing both system\footnote{As in the LTI case the appearance of $A(k)$ in \eqref{eq:noisy R(k) definition} can be interpreted as a measure of the direction of the measurement noise, which contributes to the loss of information caused \cite{DePersis2020}.} and noise information, satisfies a quadratic bound, controllers with trajectory boundedness and performance guarantees can be designed via data-dependent convex programmes, as detailed in the following subsections.

\subsection{Bounded closed-loop trajectories} 
\label{subsec:robust stab}
To ensure boundedness of the trajectories of the closed-loop system \eqref{eq:noisy closed-loop}
we derive a bound (alternative to \eqref{eq:traj_bound}), which is related to the notion of input-to-state stability. 
\begin{lemma}
Suppose there exists $P(k) = P(k)^\top$ satisfying \eqref{eq:P boundedness}, for $k = 0,\ldots, T$, and \eqref{eq:stability condition} for some $K(k)$, for $k = 0,\ldots, T-1$.
The state trajectories of the system \eqref{eq:noisy closed-loop} satisfy 
\begin{multline}
\| x(k) \| \leq \sqrt{\frac{\rho}{\eta}} \left(1-\frac{1}{\rho} \right)^{\frac{k}{2}} \| x(0) \| \\
+ \gamma_1 \left( {|v|}_{k-1} ,k \right) + \gamma_2 \left( {|d|}_{k-1},k \right),
\label{eq:noisy bound}
\end{multline}
for $k=0,\ldots, T$, with $\gamma_1(\cdot,k)$,
$\gamma_2(\cdot, k)$
class $\mathcal{K}$-functions.
\label{le:noisy closed-loop bound}
\end{lemma}
\begin{proof}
The state response at time $k$ is given by $x(k) = S_t(k,0) x(0) 
+ \sum_{j=0}^{k-1} S_t(k-1,j) ( B(j)K(j)v(j) + d(j))$,
where $S_t(k,0)$ is the state transition matrix corresponding to \eqref{eq:closed loop} as defined in Section \ref{subsec:bound}. From Theorem~\ref{th: param boundedness} we know that if there exist $P(k)=P(k)^\top$, $K(k)$ satisfying \eqref{eq:P boundedness}, for $k=0,\ldots, T$, and \eqref{eq:stability condition}, for $k = 0,\ldots, T-1$,  then 
$
\| S_t(k,0) \| \leq \sqrt{\frac{\rho}{\eta}} \left(1-\frac{1}{\rho} \right)^{\frac{k}{2}}\,,
$
for $k = 0,\ldots, T$. Combined with  properties of the operator norm this gives \eqref{eq:noisy bound} with
\begin{equation}
\begin{split}
\gamma_1(|v|_{k-1}, k) &= b \left( \sum_{j=0}^{k-1} \sqrt{\frac{\rho}{\eta}} \left( 1- \frac{1}{\rho} \right)^{\frac{k-1-j}{2}} \| K(j) \| \right) |v|_{k-1}, \\
\gamma_2(|d|_{k-1}, k) &= \left( \sum_{j=0}^{k-1} \sqrt{\frac{\rho}{\eta}} \left( 1- \frac{1}{\rho} \right)^{\frac{k-1-j}{2}} \right) |d|_{k-1},
\end{split}
\end{equation}
where $b$ denotes the upper bound on the singular values of $B(k)$, \emph{i.e.} $\| B(k) \| \leq b$ for $0 \leq j \leq k-1$.
\end{proof}

\noindent With the aim of designing 
controllers such that \eqref{eq:noisy bound} holds for $k = 0, \ldots, T$ directly using noisy data, 
we combine the result of Lemma~\ref{le:noisy closed-loop bound} and the system representation  \eqref{eq:noisy param closed loop}-\eqref{eq:noisy R(k) definition}.

\begin{theorem}
Consider the system \eqref{eq:noisy system dynamics} and suppose an ensemble of input-output data is available to form the matrices \eqref{eq:Z(k)}, \eqref{eq:U(k)}, such that the rank condition \eqref{eq:noisy rank condition} holds, for $k = 0, \ldots, T-1$. Suppose $W(k)$ satisfies
\begin{equation}
\begin{bmatrix}
I_n \\
W(k)^\top
\end{bmatrix}^\top 
\begin{bmatrix}
Q_r(k) & S_r(k) \\
S_r(k)^\top & R_r(k)
\end{bmatrix}
\begin{bmatrix}
I_n \\
W(k)^\top
\end{bmatrix} \succeq 0,
\label{eq:noisy quadratic bound}
\end{equation}
where $Q_r(k) \in \mathbb{R}^{n \times n}$, $S_r(k) \in \mathbb{R}^{n \times L}$ and $R_r(k) \prec 0 \in \mathbb{R}^{L \times L}$, for $k = 0,\ldots, T-1$. Any sequences of matrices $Y(k)$, $P(k) = P(k)^\top$ satisfying 
\begin{subequations}
\begin{equation}
\begin{bmatrix}
P(k+1) - I_n - Q_r(k) & -S_r(k) & Z(k+1)Y(k) \\
-S_r(k)^\top & -R_r(k) & Y(k) \\
Y(k)^\top Z(k+1)^\top & Y(k)^\top & P(k)  
\end{bmatrix} \succ 0,
\label{eq:noisy feasibility ineq}
\end{equation}
 \begin{equation}
 Z(k)Y(k) = P(k),
 \label{eq:noisy feasibility eq}
\end{equation}
for $k=0,\ldots,T-1$, and \eqref{eq:P boundedness},
\label{eq:noisy feasibility}
\end{subequations}
for $k= 0,\ldots, T$, where $\eta \geq 1$ and $\rho > \eta$ are finite constants, are such that the trajectories of the system \eqref{eq:noisy closed-loop}, with $K(k)$ given by \eqref{eq:param control gain}, for $k = 0,\ldots, T-1$, satisfy \eqref{eq:noisy bound}, for $k = 0,\ldots, T$.
\label{th:noisy trajectory boundedness}
\end{theorem}
\begin{proof}
By Lemma~\ref{le:noisy closed-loop bound}, \eqref{eq:noisy bound} holds for the trajectories of \eqref{eq:noisy closed-loop} if there exist $P(k) = P(k)^\top$ satisfying \eqref{eq:P boundedness}, for $k = 0, \ldots, T$, and \eqref{eq:stability condition}, for $k = 0, \ldots, T-1$.
Using \eqref{eq:noisy param closed loop}-\eqref{eq:noisy R(k) definition}, letting $Y(k) := G(k) P(k)$, and via the concrete version of the full block S-procedure (see \cite{SchererLMI}, \cite{Scherer2000}), 
\eqref{eq:stability condition} is satisfied if \eqref{eq:noisy quadratic bound} holds and $P(k)$, $Y(k)$ satisfy a quadratic matrix inequality, which can be transformed into the LMI \eqref{eq:noisy feasibility ineq} by performing the matrix multiplication, applying the Schur complement and a congruence transformation with $\text{diag} \left( I_{n+L}, P(k) \right)$. The constraint \eqref{eq:noisy feasibility eq} stems from the upper row block of \eqref{eq:noisy G definition}.
The lower row block of \eqref{eq:noisy G definition} is satisfied by $K(k)$ as in \eqref{eq:param control gain}.
\end{proof}

Quantifying \eqref{eq:noisy bound} requires knowledge of $b$ and the upper bound on the norm of the noise
vectors, $|v|_{T-1}$ and $|d|_{T-1}$.
Similarly, the condition \eqref{eq:noisy quadratic bound} cannot be verified using data alone, since $W(k)$ (as defined in \eqref{eq:noisy R(k) definition}) contains information of both the unknown system dynamics matrix and the noise affecting the data samples. 
Hence, to verify \eqref{eq:noisy quadratic bound}, knowledge of an upper bound on $A(k)$, for $k = 0, \ldots, T-1$, and the matrices $V(k)$, for $k = 0, \ldots, T$, and $D(k)$, for $k = 0, \ldots, T-1$, \emph{i.e.} the ensembles of (unmeasured) samples of measurement and process noise corresponding to the measured input-output data, is required.\footnote{While we assume that 
the system dynamics and noise are unknown, for many practical applications it is expected that reasonable upper bounds on these quantities can be estimated \cite[Chapter 8]{Khalil2017}.}
The practical relevance of the result of Theorem~\ref{th:noisy trajectory boundedness} is illustrated in \cite{Scarpa2022}, which proposes a data-driven controller for planar snake robot locomotion, partly based on this result. 

\begin{remark}
In the absence of measurement noise \eqref{eq:noisy quadratic bound} becomes a bound on $D(k)$ (the ensemble of process noise samples corresponding to the measured input-output data), which is similar to the bound on the noise data introduced in \cite{Berberich2020} for LTI systems subject to process noise only.
Note that in the LTV case \eqref{eq:noisy quadratic bound} is required to hold at each time step.
\label{re:noisy compare bound to Berberich}
\end{remark}
\begin{remark}
The result of Theorem~\ref{th:noisy trajectory boundedness} requires \eqref{eq:noisy quadratic bound} to hold only for the measured data used for the representation \eqref{eq:noisy param closed loop}-\eqref{eq:noisy R(k) definition}.
Subsequently, \eqref{eq:noisy bound} is satisfied by the trajectories of \eqref{eq:noisy closed-loop} for \emph{arbitrary}, bounded noise inputs $d(k)$, $k = 0, \ldots, T-1$, and $v(k)$, $k = 0, \ldots, T$.
\end{remark}

\begin{remark}
The matrices $Q_r(k)$, $S_r(k)$ and $R_r(k)$ in \eqref{eq:noisy quadratic bound} are chosen by the user.
This makes the quadratic bound \eqref{eq:noisy quadratic bound} a flexible condition, which contains many practical bounds as special cases, e.g. a bound on the maximum singular value (see \cite{Berberich2020}) of $W(k)$, for $k = 0,\ldots,T-1$. The choice $Q_r(k) = Z(k+1)Z(k+1)^\top$, $S_r(k) = 0$ and $R_r(k) = -\gamma(k) I_L$, for some $\gamma(k) > 0 \in \mathbb{R}$, gives the signal-to-noise ratio condition
\begin{equation}
W(k)W(k)^\top \preceq \frac{1}{\gamma(k)} Z(k+1)Z(k+1)^\top,
\label{eq:noisy SNR bound}
\end{equation}
for $k = 0, \ldots, T-1$. This condition is similar (apart from being required to hold at each time step) to the condition presented in \cite[Assumption 2]{DePersis2020} for LTI systems and represents a measure of the loss of information caused by the noise.
\label{re:SNR}
\end{remark}

\begin{figure*}[!t]
\normalsize
\setcounter{mytempeqncnt}{\value{equation}}
\setcounter{equation}{38}
\begin{subequations}
\setcounter{equation}{\value{storesubeqcounter}}
\begin{equation}
\begin{bmatrix}
\mathcal{P}(k+1) - Q_r(k) & \star & \star & \star & \star \\
- \tilde{S}_p(k)^\top \bar{E}_{cl}(k)^\top & -D_{cl}(k) \tilde{S}_p(k) - \tilde{S}_p(k)^\top D_{cl}(k)^\top + \tilde{R}_p(k) & \star & \star & \star \\
- S_r(k)^\top & - \bar{M}(k) \tilde{S}_p(k) & -R_r(k) & \star & \star \\
\bar{E}_{cl}(k)^\top & D_{cl}(k)^\top & \bar{M}(k)^\top & -\tilde{Q}_p(k)^{-1} & \star \\
\left( Z(k+1) Y(k) \right)^\top & \left( C(k) \mathcal{P}(k) + D_u(k) U(k) Y(k) \right)^\top & Y(k)^\top & 0 & \mathcal{P}(k) 
\end{bmatrix} \succ 0 
\label{eq:noisy perf data-driven performance condition inequality 1}
\end{equation}
\end{subequations}
\setcounter{equation}{\value{mytempeqncnt}}
\hrulefill
\end{figure*}
\vspace{-0.15cm}
\subsection{Robust performance} \label{subsec:robust perf}
In this subsection, we consider the problem of designing controllers of the form \eqref{eq:noisy state-feeback} for the (unknown) LTV system \eqref{eq:noisy system dynamics},
such that the closed-loop system \eqref{eq:noisy closed-loop} fulfils a disturbance attenuation condition.
Consider the performance output
\begin{equation}
\begin{split}
z(k) &= C(k)x(k) + D_u(k) u(k) + D_d(k) d(k), \\
z_f &= C(N) x(N),
\end{split}
\label{eq:noisy performance output}
\end{equation}
for $k = 0,\ldots, N-1$, where $z \in \mathbb{R}^q$, $z_f \in \mathbb{R}^q$ and $C(N), \  C(k) \in \mathbb{R}^{q \times n}$, $D_u(k) \in \mathbb{R}^{q \times m}$, $D_d(k) \in \mathbb{R}^{q \times n}$ are known matrices.
This results in the closed-loop system
\begin{subequations}
\begin{align}
x(k+1) &= A_{cl}(k) x(k) + E_{cl}(k) \bar{w}(k), \label{eq:noisy perf closed-loop x} \\
z(k) &= C_{cl}(k) x(k) + D_{cl}(k) \bar{w}(k), \label{eq:noisy perf closed-loop z}\\
z_f &= C(N) x(N), \label{eq:noisy perf closed-loop final perf output}\\
\zeta(k) &= x(k) + v(k), \label{eq:noisy perf closed-loop measurement output}
\end{align}
\label{eq:noisy perf closed-loop}%
\end{subequations}
for $k = 0, \ldots, N-1$, where $A_{cl}(k) = A(k) + B(k) K(k)$, $C_{cl}(k) = C(k) + D_u(k)K(k)$,  $\bar{w}(k):= \left[
\bar{v}(k)^\top \,\,\,\, 
d(k)^\top
\right]^\top$, $\bar v(k) = K(k)v(k)$,\,\,\,\, $ E_{cl} = \left[
B(k)\,\,\,\,  I_n 
\right]$ and  $D_{cl} = \left[
D_u(k)\,\,\,\,  D_d(k)
\right]$. 
Regarding $\bar{w}(k) \in \mathbb{R}^{(m+n)}$ as the disturbance, consider the quadratic robust performance criterion
\begin{multline}
z_f^\top z_f + \sum_{k = 0}^{N-1} \begin{bmatrix}
\bar{w}(k) \\
z(k)
\end{bmatrix}^\top \begin{bmatrix}
Q_p(k) & S_p(k) \\
S_p(k)^\top & R_p(k) 
\end{bmatrix} \begin{bmatrix}
\bar{w}(k) \\
z(k)
\end{bmatrix} \\
+ \varepsilon \sum_{k = 0}^{N-1} \bar{w}(k)^\top \bar{w}(k) \leq 0,
\label{eq:noisy perf robust performance criterion}
\end{multline}
for all $\bar{w} \in \ell_2$, where $\varepsilon > 0$ and $Q_p(k) \in \mathbb{R}^{(n+m) \times (n+m)}$, $S_p(k) \in \mathbb{R}^{(n+m) \times q}$ and $R_p(k) \succeq 0 \in \mathbb{R}^{q \times q}$, for $k = 0, \ldots, N-1$. This is the finite-horizon equivalent to the
performance criterion introduced in \cite{SchererLMI}, \cite{Scherer2000} and it captures many popular robust performance measures.
For example, the choice $Q_p(k) = - \bar{\gamma}^2 I_{(m+n)}$, $S_p(k)= 0$ and $R_p(k) = I_q$, with $\bar{\gamma} > 0 \in \mathbb{R}$, for $k = 0, \ldots, N-1$, recovers the finite-horizon $H_\infty$-control problem for discrete LTV systems (see \emph{e.g.} \cite{Ichikawa2001}).
Assuming the performance index is invertible, let
\begin{equation*}
\begin{bmatrix}
\tilde{Q}_p(k) & \tilde{S}_p(k) \\
\tilde{S}_p(k)^\top & \tilde{R}_p(k)
\end{bmatrix} = \begin{bmatrix}
Q_p(k) & S_p(k) \\
S_p(k)^\top & R_p(k)
\end{bmatrix}^{-1},
\end{equation*}
and further assume $\tilde{Q}_p(k) \prec 0$. 
The following result provides a strategy to design controllers ensuring the trajectories of \eqref{eq:noisy perf closed-loop} satisfy \eqref{eq:noisy perf robust performance criterion}. For further
results regarding robust performance of LTV systems see \emph{e.g.} \cite{Ichikawa2001,fry2017iqc}.
\begin{lemma}
Suppose there exists a matrix sequence $\mathcal{P}(k) = \mathcal{P}(k)^\top \succ 0$ satisfying
\begin{subequations}
\begin{multline}
\begin{bmatrix}
\star & \star \\
\star & \star \\
\cmidrule(lr){1-2}
\star & \star \\
\star & \star
\end{bmatrix}^\top \begin{bmatrix}[cc|cc]
- \mathcal{P}(k) & 0 & 0 & 0 \\
 0 & \mathcal{P}(k+1) & 0 & 0 \\
 \cmidrule(lr){1-4}
 0 & 0 & \tilde{Q}_p(k) & - \tilde{S}_p(k) \\
 0 & 0 & - \tilde{S}_p(k)^\top & \tilde{R}_p(k)
\end{bmatrix} \\
\times \begin{bmatrix}
A_{cl}(k)^\top & C_{cl}(k)^\top \\
I_n & 0 \\
\cmidrule(lr){1-2}
E_{cl}(k)^\top & D_{cl}(k)^\top \\
0 & I_p
\end{bmatrix} \succ 0,
\label{eq:noisy performance condition k}
\end{multline}
for $k = 0, \ldots, N-1$, and
\begin{equation}
I_p - C(N) \mathcal{P}(N) C(N)^\top \succeq 0.
\label{eq:noisy performance condition final}
\end{equation} 
\label{eq:noisy performance condition} %
\end{subequations}
The output $z(k)$ of the closed-loop system \eqref{eq:noisy perf closed-loop} subject to the disturbance input $\bar{w}(k)$ and with initial condition $x(0) = 0$ satisfies the quadratic robust performance criterion \eqref{eq:noisy perf robust performance criterion}.
\label{le:robust performance}
\end{lemma}
\begin{proof}
The result follows from dissipativity arguments (see \emph{e.g.} \cite{Byrnes1994}) and the dualisation lemma \cite[Lemma 4.9]{SchererLMI}. 
\end{proof}
With the aim of designing controllers, such that \eqref{eq:noisy perf robust performance criterion} holds directly using noisy data, consider \eqref{eq:noisy param closed loop}-\eqref{eq:noisy R(k) definition}. 
A complication then arises due to the fact that we consider measurement noise in addition to process noise.
Namely, $E_{cl}(k)$, through which the disturbance input $\bar{w}(k)$ enters the system \eqref{eq:noisy perf closed-loop} depends on the unknown input matrix $B(k)$. Hence, \eqref{eq:noisy perf closed-loop} cannot be represented using \eqref{eq:noisy param closed loop}-\eqref{eq:noisy R(k) definition} alone. To address this, we introduce an 
additional data-driven representation of $B(k)$. Supposing \eqref{eq:noisy rank condition} holds, 
$B(k)$ can be written as
\begin{equation}
B(k) = \begin{bmatrix}
A(k) & B(k)
\end{bmatrix} \begin{bmatrix}
0 \\
I_m
\end{bmatrix} = \left( Z(k+1) + W(k) \right) M(k),
\label{eq:noisy B(k) parametrisation}
\end{equation} 
with $M(k) \in \mathbb{R}^{L \times m}$ satisfying
\begin{equation}
\begin{bmatrix}
0 \\
I_m 
\end{bmatrix} = \begin{bmatrix}
Z(k) \\
U(k)
\end{bmatrix} M(k),
\label{eq:noisy M(k) definition}
\end{equation}
for $k = 0, \ldots, N-1$.
Using \eqref{eq:noisy param closed loop}-\eqref{eq:noisy R(k) definition} and \eqref{eq:noisy B(k) parametrisation}-\eqref{eq:noisy M(k) definition} the
system \eqref{eq:noisy perf closed-loop} can be equivalently written as a data-dependent lower linear fractional transformation (LFT, see \emph{e.g.} \cite{Green_Michael2012}), namely
\begin{equation}
\begin{split}
\begin{bmatrix}
x(k+1) \\
z(k) \\
\tilde{z}(k)
\end{bmatrix} &= \begin{bmatrix}
Z(k+1)G(k) & \bar{E}_{cl}(k) & I_n \\
C_{cl}(k) & D_{cl}(k) & 0 \\
G(k) & \bar{M}(k) & 0 
\end{bmatrix} \begin{bmatrix}
x(k) \\
\bar{w}(k) \\
\tilde{w}(k)
\end{bmatrix},
\end{split}
\label{eq:noisy perf LFT}
\end{equation} 
where $\tilde{w}(k) = W(k) \tilde{z}(k)$,  $\bar{E}_{cl}(k) = \left[
Z(k+1)M(k)\,\,\,\,  I_n
\right]$ and 
$\bar{M}(k) = \left[
M(k)\,\,\,\, 0
\right]$,
together with \eqref{eq:noisy perf closed-loop final perf output}, \eqref{eq:noisy perf closed-loop measurement output}, for $k  = 0, \ldots, N-1$. Using this data-dependent system representation and the result of Lemma~\ref{le:robust performance}, controllers ensuring the criterion \eqref{eq:noisy perf robust performance criterion} holds
can be designed directly using noisy data. 
\begin{theorem}
Consider the system \eqref{eq:noisy system dynamics} and suppose an ensemble of input-output data is available to form the matrices \eqref{eq:Z(k)}, \eqref{eq:U(k)}, such that the rank condition \eqref{eq:noisy rank condition} holds, for $k = 0, \ldots, N-1$. Suppose $W(k)$, as defined in \eqref{eq:noisy R(k) definition}, satisfies \eqref{eq:noisy quadratic bound}, for $k = 0, \ldots, N-1$.
Any sequences of matrices $Y(k)$, $\mathcal{P}(k) = \mathcal{P}(k)^\top$ satisfying \eqref{eq:noisy perf data-driven performance condition inequality 1} (defined above),
\begin{subequations}
\addtocounter{equation}{1}%
\setcounter{storesubeqcounter}{\value{equation}}%
\begin{equation}
Z(k)Y(k) = \mathcal{P}(k),
\label{eq:noisy perf data-driven performance condition equality 1}
\end{equation}
\begin{equation}
\begin{bmatrix}
0 \\
I_m
\end{bmatrix} = \begin{bmatrix}
Z(k) \\
U(k)
\end{bmatrix} M(k),
\label{eq:noisy perf data-driven performance condition equality 2}
\end{equation}
for $k = 0, \ldots, N-1$, and \eqref{eq:noisy performance condition final}, are such that the trajectories of the system \eqref{eq:noisy perf closed-loop}, with
\label{eq:noisy perf data-driven performance condition}
\end{subequations}
\begin{equation}
K(k) = U(k) Y(k) \mathcal{P}(k)^{-1},
\label{eq:noiy perf data control law} 
\end{equation}
for $k = 0, \ldots, N-1$, and with initial condition $x(0) = 0$, satisfy the quadratic robust performance criterion \eqref{eq:noisy perf robust performance criterion}.
\label{th:robust performance}
\end{theorem}
\begin{proof}
By Lemma~\ref{le:robust performance}, \eqref{eq:noisy perf robust performance criterion} is satisfied for trajectories of \eqref{eq:noisy perf closed-loop} if there exists $\mathcal{P}(k) = \mathcal{P}(k)^\top \succ 0$ such that \eqref{eq:noisy performance condition} holds, for $k = 0,\ldots, N-1$. 
Consider the system representation \eqref{eq:noisy perf LFT} (based on \eqref{eq:noisy param closed loop}-\eqref{eq:noisy R(k) definition} and \eqref{eq:noisy B(k) parametrisation}-\eqref{eq:noisy M(k) definition}) and let $Y(k):= G(k) \mathcal{P}(k)$. 
Via the concrete version of the full block S-procedure (\cite{SchererLMI}, \cite{Scherer2000}), \eqref{eq:noisy performance condition k} is satisfied if \eqref{eq:noisy quadratic bound} holds and $\mathcal{P}(k)$, $Y(k)$ satisfy a quadratic matrix inequality, which can be transformed into the LMI \eqref{eq:noisy perf data-driven performance condition inequality 1} by performing the matrix multiplication and applying the Schur complement twice. The equality constraints \eqref{eq:noisy perf data-driven performance condition equality 1} and \eqref{eq:noisy perf data-driven performance condition equality 2} stem from the upper row block of \eqref{eq:noisy G definition} and \eqref{eq:noisy M(k) definition}, respectively, while the lower row block of \eqref{eq:noisy G definition} is automatically satisfied by $K(k)$ in
\eqref{eq:noiy perf data control law}. 
\end{proof}

Theorem~\ref{th:robust performance} provides a general approach to design controllers guaranteeing robust quadratic performance for unknown LTV systems, affected by both measurement and process noise, directly using noisy data.

\begin{remark}
While \eqref{eq:noisy B(k) parametrisation}-\eqref{eq:noisy M(k) definition} correspond to uniquely identifying the matrix $B(k)$,  
since the sequence $M(k)$ is determined \emph{at the same time} as the control gain $K(k)$ in \eqref{eq:noiy perf data control law}, the result of Theorem~\ref{th:robust performance} is still a direct data-driven control approach (as opposed to indirect approaches involving \emph{sequential} system identification and control design).
\end{remark}

\begin{remark}
In the absence of measurement noise, \eqref{eq:noisy quadratic bound} reduces to a bound on $D(k)$ (as discussed in Remark~\ref{re:noisy compare bound to Berberich}). Moreover, the closed-loop system is described by \eqref{eq:noisy perf closed-loop x}-\eqref{eq:noisy perf closed-loop final perf output} with the disturbance defined as $\bar{w}(k) := d(k)$ and hence $E_{cl}(k) = I_n$ and $D_{cl}(k) = D_d(k)$. This removes the need to represent $B(k)$ via \eqref{eq:noisy B(k) parametrisation}-\eqref{eq:noisy M(k) definition}. 
The closed-loop system can be represented directly using \eqref{eq:noisy param closed loop}-\eqref{eq:noisy R(k) definition} via the LFT \eqref{eq:noisy perf LFT} with $\bar{E}_{cl}(k) = I_n$ and $\bar{M}(k) = 0$ and the data-dependent feasibility problem in Theorem~\ref{th:robust performance} reduces to finding sequences of matrices $Y(k)$ and $\mathcal{P}(k) = \mathcal{P}(k)^\top$ satisfying \eqref{eq:noisy perf data-driven performance condition inequality 1}-\eqref{eq:noisy perf data-driven performance condition equality 1}, for $k = 0,\ldots,N-1$, and \eqref{eq:noisy performance condition final}. The control law guaranteeing \eqref{eq:noisy perf robust performance criterion} for the system \eqref{eq:noisy perf closed-loop} is given by $u(k) = K(k)x(k)$, with $K(k)$ given by \eqref{eq:noiy perf data control law}.
\label{re:noisy perf process noise only}
\end{remark}
\vspace{-0.65cm}
\section{Periodically time-varying systems}\label{sec:periodic}
Consider the special case in which the time-variation of the matrices $A(k)$ and $B(k)$ is $\phi$-periodic, for some $\phi \in \mathbb{N}$, \emph{i.e.} 
\begin{equation}
\begin{aligned}
A(k+ \phi) = A(k), \
B(k+ \phi) = B(k), 
\label{eq:periodic system}
\end{aligned}
\end{equation}
for all $k \geq 0$. While the system matrices are assumed to be unknown, the periodic nature and period $\phi$ of the system may be known a priori. 
Exploiting periodicity, the requirement for an ensemble of $L$ data sequences (Assumption \ref{as:conditions}) can then be replaced by the requirement of one sufficiently long input-state data sequence capturing $L$ periods, \emph{i.e.} covering the time interval $k = 0,\ldots,\phi L$.
Moreover, this data sequence can be used to derive a data-driven system representation \emph{beyond} the interval $k = 0,\ldots,\phi L$. These observations allow us to derive \emph{infinite-horizon} results based on \emph{finite-horizon} data. 
Thus, in the following we consider the infinite-horizon versions of the control problems considered in Section~\ref{sec:no_noise} and Section~\ref{sec:H_inf} in the context of periodically time-varying systems.
\vspace{-0.35cm}

\subsection{Stabilisation} 
Exploiting periodicity, stabilising controllers for linear periodically time-varying systems 
can be designed using only a \emph{single, finite-length} data sequence.

\begin{corollary}  
Consider the linear periodically time-varying system \eqref{eq:systemdynmaics}, \eqref{eq:periodic system} and suppose input-state data is available to form the matrices \eqref{eq:X(k)}, \eqref{eq:U(k)}\footnote{Such data may stem from a single experiment
of length $L\phi$, or from an ensemble of $L$ experiments of length $\phi$.}, such that the rank condition \eqref{eq:rank condition} holds, for $k = 0,\ldots, \phi-1$. Any matrix sequences $Y(k)$, $P(k) = P(k)^\top$ satisfying \eqref{eq:param stability condition}, for $k = 0, \ldots, \phi-1$, where $\eta \geq 1$ and $\rho > \eta$ are finite constants, and \begin{equation}
	P(\phi) = P(0), \label{eq:period P periodic}
\end{equation}  
are such that the system \eqref{eq:closed loop}, \eqref{eq:periodic system}, with $K(k)$ given by \eqref{eq:param control gain}, for $k = 0, \ldots, \phi - 1$, and $K(k+n_p\phi) = K(k)$, for all $n_p \geq 0$, is exponentially stable.
\label{co:period param system stability}
\end{corollary}
\begin{proof}
The closed-loop LTV system \eqref{eq:closed loop} is exponentially stable if and only if there exists $P(k) = P(k)^\top$ satisfying \eqref{eq:param stability inequality} and \eqref{eq:stability condition} for some $K(k)$ for all $k \geq 0$.
If the system dynamics are $\phi$-periodic
the system is exponentially stable if and only if there exists a $\phi$-periodic solution $P(k)$, $K(k)$ to \eqref{eq:stability condition} \cite[Chapter 3.1]{Ichikawa2001}. Hence, we only need to find $K(k)$,
$P(k)$ satisfying \eqref{eq:stability condition} for one period, \emph{i.e.} for $k = 0, \ldots, \phi$.
Using \eqref{eq: param closed loop}, \eqref{eq: G definition} and following steps similar as in the proof of Theorem~\ref{th: param boundedness}, \eqref{eq:stability condition} is equivalent to~\eqref{eq:param stability condition}, with the additional constraint \eqref{eq:period P periodic} in place to ensure that $P(k)$ is periodic.
\end{proof}

\begin{remark}\label{re:noisy periodic system stability}
Corollary~\ref{co:period param system stability} is the infinite-horizon equivalent of Theorem~\ref{th: param boundedness} for periodically time-varying systems. Similarly, 
(noise) input-to-state stabilising controllers can be designed using noisy data by solving \eqref{eq:noisy feasibility}, \eqref{eq:P boundedness}, for $k=0,\ldots,\phi-1$, with the additional constraint \eqref{eq:period P periodic}, supposing \eqref{eq:noisy quadratic bound} holds.
This represents the infinite-horizon counterpart to Theorem~\ref{th:noisy trajectory boundedness}.
\end{remark}
\vspace{-0.5cm}

\subsection{Optimal Control}
\label{subsec:LQR periodic}
Consider the system \eqref{eq:systemdynmaics}, and suppose we are interested in finding a stabilising $u^*(k)$, for all $k \geq 0$, minimising
\begin{equation}
J \left( x(0), u(\cdot) \right) = \sum_{k=0}^{\infty} \left( x(k)^\top Q(k) x(k) + u(k)^\top R(k) u(k) \right),
\label{eq: inf LQR cost}
\end{equation}
with $Q(k)~=~Q(k)^\top~\succeq~0$ and $R(k)~=~R(k)^\top~\succ~0$, for all $k \geq 0$.
If \eqref{eq:periodic system} holds and $Q(k+\phi)~=~Q(k)$ and $R(k+\phi)~=~R(k)$, then the sequence of state feedback gains $K^*(k)$, $k \geq 0$, corresponding to the solution $u^*(k)$ is also $\phi$-periodic, \emph{i.e.} $K^*(k+\phi)=K^*(k)$ \cite[Chapter 3.1]{Ichikawa2001}.
Similarly to the finite-horizon case considered in Section~\ref{sec:LQR}, the described \emph{infinite-horizon} LQR problem can be formulated and solved via a convex programme involving LMI constraints \cite{Gattami2010}.
Exploiting periodicity, this can be solved directly using a \emph{single, finite-length} input-state data sequence.

\begin{corollary}
Consider the linear periodically time-varying system \eqref{eq:systemdynmaics}, \eqref{eq:periodic system} and suppose input-state data is available to form the matrices \eqref{eq:X(k)}, \eqref{eq:U(k)}, such that the rank condition \eqref{eq:rank condition} holds, for $k = 0,\ldots, \phi-1$. 
Consider the cost function \eqref{eq: inf LQR cost} with $Q(k+\phi) = Q(k)$ and $R(k+\phi) = R(k)$, for all $k \geq 0$. 
The optimal state feedback control gain sequence solving the infinite-horizon LQR problem with $u^\star(k) = K^*(k)x(k)$ is given by \eqref{eq:param gain}, for $k = 0, \ldots, \phi-1$, and $K^*(k+n_p \phi) = K^*(k)$, for all $n_p \geq 0$, with $H^*(k)$ and $S(k)^*$ the solution of
\begin{mini}[2]
    	  {\mathcal{S,H,O}}{ \sum_{k=0}^{\phi-1} \bigg( \Tr \left( Q(k) S(k) \right) + \Tr \left(O(k) \right) 			\bigg)}{\label{eq:data-dependent inf LQR problem periodic}}{}
    	  \addConstraint{\eqref{eq:data-dependent LQR problem inequality 1} - \eqref{eq:data-dependent LQR problem equality},}{}{}
    	  \addConstraint{S(\phi)}{=S(0),}
\end{mini}
    for $k = 0,\ldots, \phi-1$, where $ \mathcal{S} = \left\{ S(1), \ldots, S(\phi) \right\}$,
            $\mathcal{H} = \left\{ H(0), \ldots, H(\phi-1) \right\} $ and $
            \mathcal{O} = \left\{ O(0), \ldots, O(\phi-1) \right\}. $
    \label{co:inf LQR periodic}
\end{corollary}
\begin{proof}
The \emph{infinite-horizon} LQR problem can be recast as a convex programme (see \cite{Gattami2010}). Then, exploiting that the solution is a state feedback law and introducing \eqref{eq: param closed loop}, \eqref{eq: G definition} yields
\eqref{eq:data-dependent LQR problem} with $Q_f=0$ and $N \rightarrow \infty$,
where $H(k) := G(k)S(k)$. Recall that $K^*(k)$ for the considered problem is $\phi$-periodic \cite[Chapter 3.1]{Ichikawa2001}. 
It remains to be shown that this $\phi$-periodic solution can be recovered by solving \eqref{eq:data-dependent inf LQR problem periodic} over one period, with the additional constraint $S(\phi)=S(0)$. 
Since $K^*(k)$ is stabilising by construction, there exists a $\phi$-periodic solution $S^*(k+\phi) = S^*(k)$ satisfying \eqref{eq:data-dependent LQR problem inequality 1} and \eqref{eq:data-dependent LQR problem inequality 2} (this can be shown using analogous arguments as in the proof of Corollary~\ref{co:period param system stability}). Thus, the solution of the slack variable $O^*(k) = R(k)^\frac{1}{2} K^*(k) S^*(k) K^*(k)^\top R(k)^\frac{1}{2}$ is also $\phi$-periodic.
Hence, the constraints \eqref{eq:data-dependent LQR problem inequality 1} - \eqref{eq:data-dependent LQR problem equality} are satisfied at time $k + n_p \phi$, for all $n_p \geq 0$, if they are satisfied at time $k$.
Similarly, the optimal stage cost $I_c^*(k) = \Tr \left( Q(k) S^*(k) \right) + \Tr \left( O^*(k) \right),$ satisfies $I_c^*(k+n_p \phi) = I_c^*(k)$, for all $n_p \geq 0$. Hence, the optimal cost is given by $\sum_{k=0}^{\infty} I_c^*(k) = \lim_{n_p \rightarrow \infty} n_p \sum_{k=0}^{\phi-1} I_c^*(k)$.
Note that $\sum_{k=0}^{\phi-1} I_c^*(k)$ is the optimal cost obtained by solving \eqref{eq:data-dependent inf LQR problem periodic}. Hence, the periodic solution to the infinite-horizon LQR problem is given by $K^*(k)$, for $k = 0,\ldots,\phi - 1$, solving \eqref{eq:data-dependent inf LQR problem periodic}, and $K^*(k+n_p \phi) = K^*(k)$, for all $n_p \geq 0$.
\end{proof}
\vspace{-0.55cm}

\subsection{Robust performance}
Consider the problem of designing stabilising controllers of the form \eqref{eq:noisy state-feeback}, such that the closed-loop system \eqref{eq:noisy perf closed-loop x}-\eqref{eq:noisy perf closed-loop z}, \eqref{eq:noisy perf closed-loop measurement output} satisfies the infinite-horizon performance criterion
\begin{multline}
\sum_{k=0}^{\infty} \begin{bmatrix}
\bar{w}(k) \\
z(k)
\end{bmatrix}^\top \begin{bmatrix}
Q_p(k) & S_p(k) \\
S_p(k)^\top & R_p(k)
\end{bmatrix} \begin{bmatrix}
\bar{w}(k) \\
z(k)
\end{bmatrix} \\
+ \varepsilon \sum_{k=0}^{\infty} \bar{w}(k)^\top \bar{w}(k) \leq 0,
\label{eq:noisy infinite horizon robust perf criterion}
\end{multline}
for all $\bar{w} \in \ell_2$, with $\varepsilon > 0$, $R_p(k) \succeq 0$ and such that $\tilde{Q}_p(k) \prec 0$, for all $k \geq 0$. Suppose the system dynamics and the performance index are $\phi$-periodic, \emph{i.e.} \eqref{eq:periodic system} holds and $C(k + \phi) = C(k)$, $D_u(k + \phi) = D_u(k)$, $D_d(k+\phi) = D_d(k)$, $Q_p(k+\phi) = Q_p(k)$, $S_p(k+\phi) = S_p(k)$, $R_p(k+\phi) = R_p(k)$.
This problem can be solved via a data-driven convex programme using a \emph{single, finite-length} data sequence.

\begin{corollary}
Consider the linear periodically time-varying system \eqref{eq:noisy system dynamics}, \eqref{eq:periodic system} and suppose input-output data is available to form the matrices \eqref{eq:Z(k)}, \eqref{eq:U(k)}, such that the rank condition \eqref{eq:noisy rank condition} holds, for $k = 0,\ldots, \phi-1$. Suppose the performance index is $\phi$-periodic and $W(k)$, as defined in \eqref{eq:noisy R(k) definition}, satisfies~\eqref{eq:noisy quadratic bound}, for $k = 0, \ldots, \phi-1$.
Any sequences of matrices $Y(k)$, $\mathcal{P}(k) = \mathcal{P}(k)^\top$ satisfying \eqref{eq:noisy perf data-driven performance condition inequality 1}-\eqref{eq:noisy perf data-driven performance condition equality 2}, for $k = 0,\ldots,\phi-1$, and
\begin{equation}
\mathcal{P}(\phi) = \mathcal{P}(0),
\label{eq:noisy perf enforce periodicity}
\end{equation}
are such that the trajectories of the system \eqref{eq:noisy perf closed-loop x}, \eqref{eq:noisy perf closed-loop z}, \eqref{eq:noisy perf closed-loop measurement output}, with $K(k)$ given by \eqref{eq:noiy perf data control law}, for $k = 0, \ldots, \phi-1$, and $K(k+n_p\phi) = K(k)$, for $n_p \geq 0$, and with initial condition $x(0) = 0$, satisfy the quadratic robust performance criterion \eqref{eq:noisy infinite horizon robust perf criterion}.
\label{co:noisy performance periodic}
\end{corollary}
\begin{proof}
Analogous to Lemma~\ref{le:robust performance} it can be shown via dissipativity arguments (see \emph{e.g.} \cite{Byrnes1994}) and the dualisation lemma \cite[Lemma 4.9]{SchererLMI} that \eqref{eq:noisy infinite horizon robust perf criterion} holds, if there exist $\phi$-periodic sequences $K(k)$, $\mathcal{P}(k) = \mathcal{P}(k)^\top$ satisfying \eqref{eq:noisy performance condition k} for all $k \geq 0$. Stability is implied by the upper left block of \eqref{eq:noisy performance condition k} and the assumption that $\tilde{Q}_p(k) \prec 0$ for all $k \geq 0$.
The data-driven formulation \eqref{eq:noisy perf data-driven performance condition}, \eqref{eq:noisy perf enforce periodicity} follows via analogous steps to those in the proof of Theorem~\ref{th:robust performance}, exploiting periodicity.
\end{proof}

\section{Conclusion} \label{sec:conclusion}
A model-free, data-driven representation of closed-loop LTV systems under state feedback has been employed to design feedback controllers ensuring that the resulting closed-loop trajectories satisfy certain boundedness, performance and robustness criteria via the formulation of convex feasibility/optimisation problems involving data-dependent LMIs.
Both the noise-free case and the case in which the data and the system are affected by process and measurement noise have been considered. 
Special insights have also been provided for the case of periodically time-varying systems. 

\bibliographystyle{IEEEtran}
\bibliography{biblio}            
\end{document}